\documentclass[a4paper,12pt]{article}
\title{The Weierstrass semigroups on double covers of genus two curves}
\author{
Takeshi Harui\thanks{E-mail: takeshi@cwo.zaq.ne.jp, kt13459@ns.kogakuin.ac.jp}\\
{\small Academic Support Center, Kogakuin University}\\
{\small  Hachioji, 192-0015, Japan}\\
Jiryo Komeda\thanks{E-mail: komeda@gen.kanagawa-it.ac.jp}\\
{\small Department of Mathematics, Center for Basic Education and Integrated Learning}\\
{\small Kanagawa Institute of Technology, Atsugi, 243-0292, Japan}\\
and\\
Akira Ohbuchi\thanks{E-mail: ohbuchi@tokushima-u.ac.jp}\\
{\small Department of Mathematics, Faculty of Integrated Arts and Sciences}\\
{\small Tokushima University, Tokushima, 770-8502, Japan}\\
}

\date{}

\usepackage{amssymb}
\usepackage[leqno]{amsmath}
\usepackage{latexsym}
\usepackage{theorem}
\usepackage{enumerate}
\newcounter{examplec}[section]

\newcommand{\qed}{\hfill $  \Box $}

\newcommand{\tH}{\tilde{H}}
\newcommand{\tC}{\tilde{C}}
\newcommand{\tP}{\tilde{P}}
\newcommand{\tg}{\tilde{g}}
\newcommand{\tih}{\tilde{h}}
\newcommand{\la}{\langle}
\newcommand{\ra}{\rangle}

\newcommand{\dis}{\displaystyle}

\newcommand{\NI}{\mathbb{N}_0}

\newcommand{\DC}{\mbox{of double covering type}}

\Roman{equation}
\begin{document}
\maketitle
\renewcommand{\thefootnote}{\fnsymbol{footnote}}
\footnotetext{
The second author is partially supported by
Grant-in-Aid for Scientific Research (24540057), Japan Society for the Promotion Science.
The third author is partially supported by
Grant-in-Aid for Scientific Research (24540042), Japan Society for the Promotion Science.
 }
\renewcommand{\abstractname}{}
\begin{abstract}\vskip-2mm
We show that three numerical semigroups $\la 5,6,7,8 \ra$, $\la 3,7,8 \ra$ and $\la 3,5 \ra$ are of double covering type, i.e., the Weierstrass semigroups of ramification points on double covers of curves.
Combining the result with~\cite{oli-pim} and~\cite{kom2} we can determine the Weierstrass semigroups of the ramification points on double covers of genus two curves.

\vspace{2mm} \noindent 
\end{abstract}
\section{Introduction}
\label{intro}
Let $C$ be a complete nonsingular irreducible curve over an algebraically closed field $k$ of characteristic 0, which is called a {\it
curve} in this paper.
For a point $P$ of $C$, we set
$$H(P)=\{\alpha\in \mathbb{N}_{0}|\mbox{ there exists a rational function } f \mbox{ on }C\mbox{ with } (f)_{\infty}=\alpha P\},$$
which is called the {\it Weierstrass semigroup of $P$} where $\mathbb{N}_{0}$ denotes the additive monoid of non-negative integers.
A submonoid $H$ of $\mathbb{N}_0$ is called a {\it numerical semigroup} if its complement $\mathbb{N}_0\backslash H$ is a finite set.
The cardinality of $\mathbb{N}_0\backslash H$ is called the {\it genus} of $H$, which is denoted by $g(H)$.
It is known that the Weierstrass semigroup of a point on a curve of genus $g$ is a numerical semigroup of genus $g$.
For a numerical semigroup $\tH$ we denote by $d_2(\tH)$ the set of consisting of the elements  $\tih/2$ with even $\tih\in \tH$, which becomes a numerical semigroup.
A numerical semigroup $\tH$ is said to be {\it of double covering type} if there exists a double covering $\pi:{\tilde C}\longrightarrow C$ of a curve with a ramification point $\tP$ over $P$ satisfying $H(\tP)=\tH$.
In this case we have $d_2(H(\tP))=H(P)$.

We are interested in numerical semigroups of double covering type.
Let $\tH$ be a numerical semigroup of genus $\tg$ with $d_2(\tH)=\NI$ whose genus is $0$.
Then the semigroup $\tH$ is  $\la 2,2\tg+1\ra$, where for any positive integers $a_{1},a_{2},\ldots,a_{n}$ we denote by $\langle a_{1},a_{2},\ldots,a_{n}\rangle$ the additive monoid generated by $a_{1},a_{2},\ldots,a_{n}$.
In this case $\tH$ is the Weierstrass semigroup of a ramification point $\tP$ on a double cover of the projective line which is of genus $\tg$.
Hence, $\tH$ is \DC.

Let $\tH$ be a numerical semigroup of genus $\tg$ with $d_2(\tH)=\la 2,3 \ra$ which is the only one numerical semigroup of genus $1$.
Then the semigroup $\tH$ is either $\la 3,4,5\ra$ or $\la 3,4\ra$ or $\la 4,5,6,7\ra$ or $\langle 4,6,2\tg-3\rangle$ with $\tg\geqq 4$ or $\langle 4,6,2\tg-1,2\tg+1\rangle$ with $\tg\geqq 4$.
We can show that there is a double covering of an elliptic curve with a ramification point whose Weierstrass semigroup is any semigroup in the above ones (for example, see~\cite{kom1},~\cite{kom2}).

Oliveira and Pimentel~\cite{oli-pim} studied the semigroup $\tH=\langle 6,8,10,n\rangle$ with an odd number $n\geqq 11$.
They showed that the semigroup $\tH$ is \DC.
In this case we have $d_2(\tH)=\la 3,4,5 \ra$, which is of genus $2$.
Moreover, in~\cite{kom2} we proved that any numerical semigroup $\tH$ with $d_2(\tH)=\la 3,4,5\ra$ except $\la 5,6,7,8 \ra$, $\la 3,7,8 \ra$, $\la 3,5 \ra$ and $\la 3,5,7 \ra$ is \DC.
In view of the fact that $g(\la 3,5,7 \ra)=3<2\cdot 2$ the semigroup $\la 3,5,7 \ra$ is not \DC.
There is another numerical semigroup of genus $2$, which is $\la 2,5 \ra$.
Using the result of Main Theorem in~\cite{kom-ohb} every numerical semigroup $\tH$ with $d_2(\tH)=\la 2,5 \ra$ is \DC.
In this paper we will study the remaining three numerical semigroups.
Namely we prove the following:
\vskip3mm
\noindent
{\bf Theorem 1} {\it The three numerical semigroups $\la 5,6,7,8 \ra$, $\la 3,7,8 \ra$ and $\la 3,5 \ra$ are of double covering type.}
\vskip3mm
Combining this theorem with the results in~\cite{oli-pim} and~\cite{kom2}, we have the following conclusion:
\vskip3mm
\noindent
{\bf Theorem 2} {\it Let $\tH$ be a numerical semigroup with $g(d_2(\tH))=2$.
If $\tH\not=\la 3,5,7\ra$, then it is of double covering type.} 
\section{The proof of Theorem}
To prove that the three numerical semigroups are \DC  \ we use the following remark which is stated in Theorem 2.2 of~\cite{kom3}.
\vskip3mm
\noindent
{\bf Remark 1. }Let $\tH$ be a numerical semigroup.
We set
$$n=\min\{\tih\in \tH\mid \tih\mbox{ is odd}\}\mbox{ and }g(\tH)=2g(d_2(\tH))+\frac{n-1}{2}-r$$
with some non-negative integer $r$.
Assume that $H=d_2(\tH)$ is Weierstrass.
Take a pointed curve $(C,P)$ with $H(P)=H$.
Let $Q_1,\ldots,Q_r$ be points of $C$ different from $P$ with $h^0(Q_1+\cdots+Q_r)=1$.
Moreover, assume that $\tH$ has an expression
$$\tH=2H+\la n,n+2l_1,\ldots,n+2l_s\ra$$
of generators with positive integers $l_1,\ldots,l_s$ such that
$$h^0(l_iP+Q_1+\cdots+Q_r)=h^0((l_i-1)P+Q_1+\cdots+Q_r)+1$$
for all $i$.
If the divisor $nP-2Q_1-\cdots-2Q_r$ is linearly equivalent to some reduced divisor not containing $P$, then there is a double covering $\pi:\tC\longrightarrow C$ with a ramification point $\tP$ over $P$ satisfying $H(\tP)=\tH$, hence $\tH$ is $\DC$.
\vskip3mm
\noindent
By seeing the proof of Theorem 2.2 in~\cite{kom3} we may replace the assumption in Theorem 2.2 in~\cite{kom3} that the complete linear system $|nP-2Q_1-\cdots-2Q_r|$ is base point free by the above assumption that the divisor $nP-2Q_1-\cdots-2Q_r$ is linearly equivalent to some reduced divisor not containing $P$.
\vskip3mm
\noindent
{\it Case 1.} Let $\tH=\la 5,6,7,8\ra$.
Then we have $H=d_2(\tH)=\la 3,4,5\ra$ and $\dis g(\tH)=5=2\cdot 2+\frac{5-1}{2}-1$.
Moreover, we have $\tH=2H+\la 5,5+2\cdot 1\ra$.
Let $C$ be a curve of genus $2$ and $\iota$ the hyperelliptic involution on $C$.
Let us take a point $P$ of $C$ with $H(P)=\la 3,4,5\ra$ and $3(P-\iota(P))\not\sim 0$.
Then we get $h^0(P+\iota(P))=2=h^0(\iota(P))+1$.
Moreover, we have $R\not=P$ if the complete linear system $|5P-2\iota(P)|$ has a base point $R$.
Indeed, we assume that $R=P$.
Then we have
$$h^0(5P-2\iota(P)-P)=h^0(5P-2\iota(P))=3+1-2=2,$$
which implies that
$$4P-2\iota(P)\sim g_2^1\sim P+\iota(P).$$
Hence, we get $3(P-\iota(P))\sim 0$.
This is a contradiction.

We assume that $|5P-2\iota(P)|$ has a base point $R$.
Then we get $5P-2\iota(P)\sim R+E$, where $E$ is an effective divisor of degree $2$ with projective dimension $1$.
In this case the complete linear system $|E|$ is base point free. Therefore, the divisor $5P-2\iota(P)$ is linearly equivalent to some reduced divisor not containing $P$.
If $|5P-2\iota(P)|$ is base point free, then the divisor $5P-2\iota(P)$ satisfies the above condition.
By Remark 1 the semigroup $\tH=\la 5,6,7,8\ra$ is \DC.
\vskip3mm
\noindent
{\it Case 2.} Let $\tH=\la 3,7,8\ra$.
Then we have $H=d_2(\tH)=\la 3,4,5\ra$ and $\dis g(\tH)=4=2\cdot 2+\frac{3-1}{2}-1$.
Moreover, we have $\tH=2H+\la 3,3+2\cdot 2\ra$.
Let $C$ be a curve of genus $2$ and $\iota$ the hyperelliptic involution on $C$.
We take a point $P$ of $C$ with $H(P)=\la 3,4,5\ra$.
Let $\varphi:C\longrightarrow \mathbb{P}^1$ be a covering of degree $3$ corresponding to the complete linear system $|3P|$.
We may take the pointed curve $(C,P)$ such that $\varphi$ has a simple ramification point $Q$.
Then there is another simple ramification point
of $\varphi$ by Riemann-Hurwitz formula.
Hence, we may assume that $\iota P\not=Q$, which implies that $P+Q\not\sim g_2^1$.
Thus, we get $h^0(2P+Q)=2=h^0(P+Q)+1$.
Let $R$ be the point satisfying $2Q+R\sim 3P$.
Then we have $R\not=P$ and $3P-2Q\sim R$.
By Remark 1 the semigroup $\tH=\la 3,7,8\ra$ is \DC.
\vskip3mm
\noindent
{\it Case 3.} Let $\tH=\la 3,5\ra$.
Then we have $H=d_2(\tH)=\la 3,4,5\ra$ and $\dis g(\tH)=4=2\cdot 2+\frac{3-1}{2}-1$.
Moreover, we have $\tH=2H+\la 3,3+2\cdot 1\ra$.
Let $C$ be a curve whose function field is $k(x,y)$ with an equation $y^3=(x-c_1)(x-c_2)(x-c_3)^2$, where $c_1$, $c_2$ and $c_3$ are
distinct elements of $k$.
Let $\pi:C\longrightarrow \mathbb{P}^1$ be the morphism corresponding to the inclusion $k(x)\subset k(x,y)$.
Then $C$ is of genus $2$.
Let $P=P_1$, $P_2$, $P_3$ and $P_4$ be the ramification points of $\pi$.
Since $\pi$ is a cyclic covering, it induces an automorphism $\sigma$ of $C$ with $C/\la \sigma \ra\cong \mathbb{P}^1$.
Let $\iota$ be the hyperelliptic involution on $C$.
Then we have $\sigma\circ\iota=\iota\circ\sigma$.
Indeed, we have
$$(\sigma\circ\iota\circ\sigma^{-1})\circ(\sigma\circ\iota\circ\sigma^{-1})=\sigma\circ\iota\circ\iota\circ\sigma^{-1}=\sigma\circ\sigma^{-1}=id.$$
Hence, the automorphism $\sigma\circ\iota\circ\sigma^{-1}$ is an involution.
Moreover, we have a bijective correspondence between the sets Fix$(\iota)$ and Fix$(\sigma\circ\iota\circ\sigma^{-1})$
sending $Q$ to $\sigma(Q)$, where Fix$(\iota)$ and Fix$(\sigma\circ\iota\circ\sigma^{-1})$ are the sets of the fixed points
by $\iota$ and $\sigma\circ\iota\circ\sigma^{-1}$ respectively.
Hence, $\sigma\circ\iota\circ\sigma^{-1}$ is also the hyperelliptic involution.
Thus, we have $\sigma\circ\iota\circ\sigma^{-1}=\iota$.
Since $\sigma(\iota(P))=\iota(\sigma(P))=\iota(P)$, the point $\iota(P)$ is a fixed point of $\sigma$.
Hence, we may assume that $\iota(P)=P_2$.
Then we obtain $h^0(P+P_2)=2=h^0(P_2)+1$.
Moreover, we have
$$3P-2P_2\sim 3P_2-2P_2=P_2.$$
By Remark 1 the semigroup $\la 3,5\ra$ is \DC.
\qed

\end{document}